# Supply Chain Network Design Heuristics for Capacitated Facilities under Risk of Supply Disruptions

AliReza Madadi[*], Mary E. Kurz, Scott J. Mason, Kevin M. Taaffe

Department of Industrial Engineering, Clemson University, 110 Freeman Hall, Clemson, SC 29634, USA

**Abstract**

Recent events such as the Heparin tragedy, in which patients lost their lives due to tainted pharmaceuticals, highlight the necessity for supply chain designers and planners to consider the risk of even low probability disruptions in supply chains. One of the most effective ways to hedge against supply chain network disruptions is to robustly design the supply chain network. This involves both strategic decisions (e.g., which suppliers to source from, plant locations, etc.) and tactical decisions (e.g., capacity allocation, etc.). Since disruptions are modeled as events which occur randomly and may have a random length, we consider a mixed-integer stochastic model. However, such network design problems belong to the class of NP-hard problems. Accordingly, we develop efficient heuristic algorithms and a metaheuristic approach to obtain acceptable solutions to these types of problems in reasonable runtimes so that the process of decision making becomes facilitated with no drastic sacrifices in solution quality.

## 1. Introduction

Some supply chain disruptions are not only costly, but may have catastrophic consequences in spite of their low probability of occurrence. For instance, in the healthcare supply chain, it is not acceptable to experience a late delivery or product shortage if patients' lives are in danger. Nevertheless, several examples of disruptions in healthcare supply chains exist. For example, the disruption of a flu vaccine manufacturer in Bristol, UK in 2004 resulted in disastrous consequences. The U.K. government stopped production when U.S. regulators inspected a manufacturing plant and found evidence of bacterial contamination problems. This reduced the U.S.'s supply of the vaccine by nearly 50% during the 2004-2005 flu season [5]. A healthcare supply chain is also very susceptible to disruptions caused by contamination. Heparin, a widely-used blood-thinning medicine that is made from pig intestines, was contaminated by an undetected outbreak of blue ear pig disease in China in 2008. This led to 81 patient deaths as well as hundreds of allergic reactions in the United States [7]. The investigation engaged several government agencies, university researchers and a biotech company that has a generic heparin under FDA review. Although no one at the time knew what was causing the reactions, members of Congress concluded that the issue was the result of "regulatory failure" based on news reports that FDA had not inspected a Chinese heparin production facility [1]. In another supply chain disruption, a baby food producer who purchased vitamin supplements from a Chinese supplier found out that the supplements were contaminated by cement [8]. This incident involved 22 Chinese and 10 global manufacturers and led to kidney problems and kidney stones in Chinese babies, illustrated the result of poor or failed to inspection. [9].

---

[*] Corresponding author. Tel.: +1 864 6564716; fax: +1 864 6560795.
  E-mail address: amadadi@clemson.edu (A. Madadi).

These types of incidents accentuate the need to consider supply chain risks and (supply) disruptions in the design and planning stages. However, managers, deceived by the small likelihood of a disruption, often tend to underestimate the impact of such mishaps. This is reflected in designing supply chain networks which only take operational efficiencies into account. Unfortunately, once a disruption occurs, there are few opportunities to change existing supply chain infrastructure. Therefore, to hedge against supply chain disruptions, it is critical to consider potential disruptions during the design of the supply chain networks so that the network can be responsive and resilient in the event of an unplanned disruption.

The goal of this research is to design a single-period, single-product supply chain network model with capacitated facilities to hedge against the risk of sending tainted materials to consumers. We focus on supply disruptions causing the loss of all or a significant fraction of the production at a set of facilities in the same geographic.

Typically, the decision making process dealing with supply chain disruptions involves both strategic and tactical considerations [2]. Strategic decisions comprise decisions such as choosing which markets to serve, from which suppliers to source, the location of facilities/suppliers, how many suppliers to use, etc. Tactical decisions include decisions such as capacity allocation, inventory management and transportation planning [2].

Our model consists of facility selection and capacity allocation among facilities. A key parameter in our model is the consideration of facility inspection. This aspect of the work was inspired by tragedies such as the Heparin and Chinese baby food manufacturing incidents. If the risk of shipping tainted materials can be minimized prior to such tragedies, producers can decrease liability and improve consumer safety. Insights into how our model should be configured to avoid the risk of tainted products reaching consumers are of interest to several types of supply chains such as healthcare, pharmaceutical, cosmetic and beauty, and food or dairy industries.

The objective of the model is to minimize the expected overall cost which is composed of the cost of selecting the facility, shipping untainted products, shipping tainted products, inspect the facility, and discarding tainted products. We can formulate this problem using a two-stage stochastic mixed-integer problem. However, while this approach may allow for exact solutions in some situations, it can be very challenging to draw concrete analytical insights from such models and to obtain good solutions for large instances within a limited time frame since the problem is a special case the two-stage stochastic capacitated facility location problem which is well known to be NP-hard [3,4]. Based on our experience in solving various size problems using commercial software in this paper, we show that the number of facilities used have a significant impact on the solution time. As a result, we develop several heuristics and metaheuristic to efficiently solve and handle large size problems.

The paper is organized as follows. In Section 2, the problem description is discussed. The mathematical formulation is introduced in Section 3. In Sections 4 and 5, data generation method and solution procedure are presented, respectively. Computational results are discussed in Section 6. Finally, Section 7 presents brief conclusions.

**2.     Problem Description and Mathematical Formulation**

2.1.    Description

The earliest work in supply chain network design was developed by Geoffrion and Graves [5]. They introduced a multi-commodity logistics network design model for optimizing finished product flows from plants to distribution centers to final consumers. Beginning with the work of Geoffrion and Graves, a large number of optimization-based approaches have been proposed for the design of supply chain networks. These works have resulted in significant improvements in the modeling of these problems as well as in algorithmic and computational efficiency. However, generally this research assumes that the design parameters for the supply chain are deterministic [2,6,7,8,9,10]. Unfortunately, critical parameters such as consumers' demand, resource supply, and price of the material are generally uncertain. Therefore, traditional deterministic optimization is not suitable for truly capturing the behavior of the real-world problem.

The significance of uncertainty has encouraged a number of researchers to address stochastic parameters in their research. However, most of the stochastic approaches for supply chain network design only consider tactical level decisions usually related to demand uncertainty [11,12,13,14], while supply uncertainty is often ignored and supply capacity assumed to be unlimited. In contrast with the most prior research, we focus on the supply (capacity) management required to mitigate the impact of facility capacity disruptions. Moreover, we assume that supply quantities can be influenced by inspection, which might be conducted at facility locations.

We utilize a mixed-integer stochastic programming model that is formulated as a two-stage optimization problem. The selection of the facilities is considered at the first stage and modeled as a binary decision. The second-stage decision variables include tactical decisions which are made after realization of the random events (i.e., supply disruption). The second-stage decisions indicate the actual capacity allocation as well as inspection decisions at each facility. The inspection decision is modeled as a binary variable for each facility. Therefore, the model enables us to determine when and where inspections should be performed with the intent of minimizing the amount of tainted product shipped to consumers.

In Figure 1 we provide a hypothetical supply network with an initial assignment of consumers to facilities at a point in time before any disruptions has occurred. We consider a set of facilities and consumers. In some cases disruptions can be a consequence of tainted raw material (received from suppliers.) Hence, for sake of clarity and in order to show the flow of raw materials from the suppliers to the facility, we have illustrated the set of suppliers as well. Three facilities were selected and the capacity was sufficient to fulfill all the demand of all the consumers.

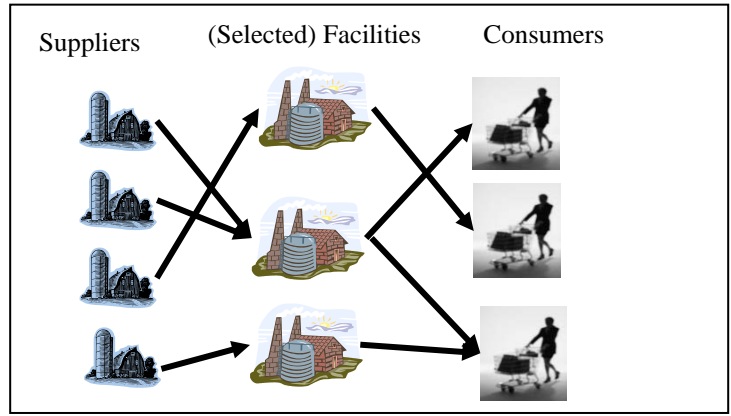

Figure 1 Initial demand allocation (before disruption)

In Figure 2, a scenario with disruptions at two facilities is presented. This disruption caused the facilities to produce tainted items and ship them to consumers.

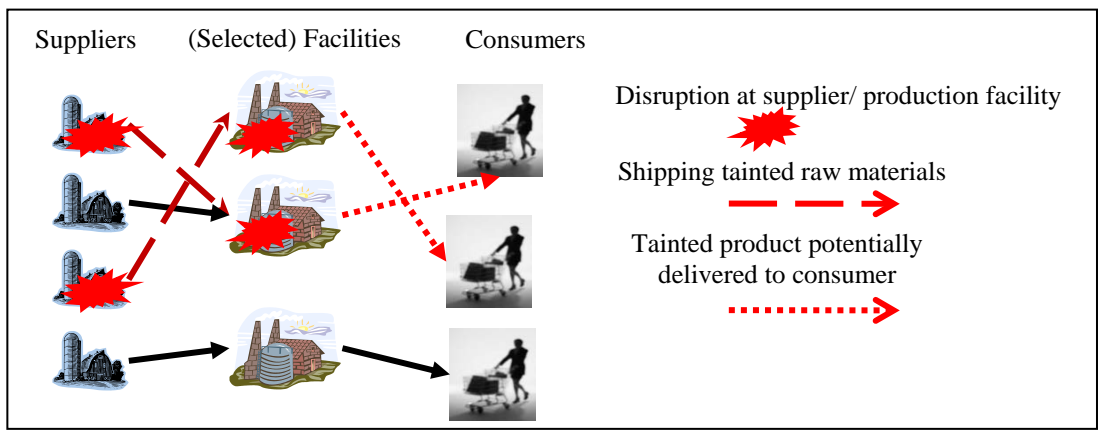

Figure 2 Shipping tainted products to consumers after disruptions (no inspection)

Once inspection is implemented in a facility, a portion of tainted items is discarded and fewer tainted items are delivered to the consumers. However, discarding the tainted items might result in consumer demand being unsatisfied. Then, the unmet demand can be fulfilled by adding another facility, as illustrated in Figure 3.

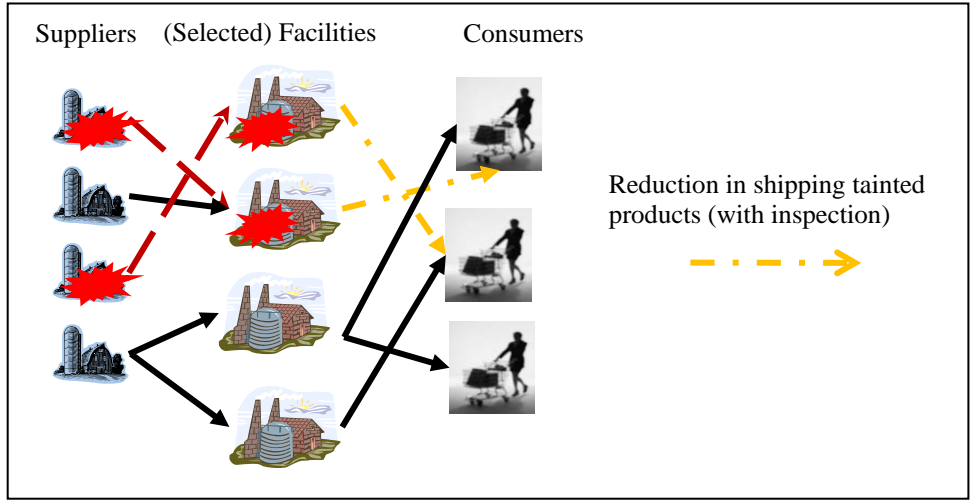

Figure 3 Reduction of delivered tainted products (after inspection)

## 2.2. Mathematical Model

Consider a supply chain network $\mathcal{N} = (L, C)$, where $L$ is the set of facilities and $C$ is the set of consumers. In the first stage, $x_l$ is 1 if facility $l$ is selected and 0 otherwise (where $l \in L$ is an index for facilities). Let $Q(x, \tilde{s})$ represents the optimal solution of the second stage problem corresponding to the first stage decision variable $x$, and the random scenario $\tilde{s}$. The stochastic formulation of the problem can be written as follows:

$$\min \sum_{l \in L} x_l f_l + \mathrm{E}\left[Q(x, \tilde{s})\right]. \tag{1}$$

$$\text{s.t. } x_l \in \{0,1\}, \forall l \in L, \tag{2}$$

where $\mathrm{E}\left[Q(x, \tilde{s})\right]$ is the expected cost taken with respect to random scenario $\tilde{s}$. The objective (1) in the first-stage problem is the sum of the cost of selecting facilities. The first-stage constraint (2) restricts the decision variables $x_l$ to be binary. Given a feasible first-stage solution vector $x$, the objective of the second-stage problem for random scenario $\tilde{s}$ minimizes the sum of the allocation (shipping) cost of the untainted products, the cost of shipping tainted product, the cost of discarding tainted product after inspection, and the cost of inspection. In this model, we discard tainted products. An alternative is to repair (which is considered as reworking) the tainted product which we may consider in future research.

To deal with the uncertainty in the second stage, a scenario-based modeling approach is proposed which has been used in stochastic programing problems [4,12]. In the second stage, let's consider random scenario $\tilde{s}$ to have a discrete distribution $\Pr(\tilde{s} = s) = \rho_s$ where $\rho_s$ is the probability of occurrence for scenario $s$. Given a finite set of scenarios, with associated probabilities $\rho_s$, $\mathrm{E}[Q(x,\tilde{s})]$ can be evaluated as $\mathrm{E}[Q(x,\tilde{s})] = \sum_{s \in S} \rho_s Q(x,s)$. Hence, we can present the deterministic equivalent of the formulation (1). To simplify, we denote this as the Supply Chain Design (SCD) model. We first summarize the complete notation for the SCD in the following.

**Sets**

| | |
|---|---|
| $C$ | the set of consumers, indexed by $c$ |
| $L$ | the set of candidate facility, indexed by $l$ |
| $S$ | the set of realized scenarios, indexed by $s$ |

**Parameters**

| | |
|---|---|
| $f_l$ | fixed cost of selecting facility $l$ |
| $\kappa_l$ | capacity of facility $l$ |
| $b_c$ | total demand of consumer $c$ |
| $n_l$ | fixed cost of implementing an inspection at candidate location $l$ |
| $q_{ls}$ | fraction of tainted products produced at facility $l$ in scenario $s$ |
| $r_{ls}$ | fraction of tainted products produced at facility $l$ after inspection (we assume $q_{ls} > r_{ls}$) in scenario $s$ |
| $\lambda_{lc}$ | cost of shipping an untainted product from facility $l$ to consumer $c$ |
| $o_{lc}$ | penalty cost for shipping a tainted product from facility $l$ to consumer $c$ |
| $\gamma_{lc}$ | cost of discarding a tainted product at facility $l$ after inspection originally destined for consumer $c$ |
| $\rho_s$ | probability of occurrence for scenario $s$ |

To make the definitions of $q_{ls}$ and $r_{ls}$ clearer, consider the following example. Suppose that under scenario $s$, the extent of failures at unreliable facility $l$ is given by $q_{ls} = 0.20$ and $r_{ls} = 0.05$. This means that for every 100 items of item produced at facility $l$, $100 q_{ls} = 20$ of them will be tainted. If no inspection is used, these 20 tainted items might be shipped to consumers. If inspection is used, 15 of these 20 tainted items will be detected and discarded while $100 r_{ls} = 5$ of them will go undetected and be shipped to consumers.

**Decision Variables**

$$x_l = \begin{cases} 1, & \text{if facility } l \text{ is selected,} \\ 0, & \text{otherwise} \end{cases}$$

$$z_l = \begin{cases} 1, & \text{if inspection is used at facility } l, \\ 0, & \textit{otherwise} \end{cases}$$

$p_{lcs}$     number of untainted products shipped from facility $l$ to consumer $c$ in scenario $s$

$k_{lcs}$     number of tainted products produced at facility $l$ intended to be shipped to consumer $c$ in scenario $s$

$d_{lcs}$     number of tainted products discarded after inspection in scenario $s$

We now present the deterministic equivalent of the formulation. Note that the second-stage decision variables are indexed by a scenario index.

[SCD]    $\min \sum_{l \in L} x_l f_l + \sum_{s \in S} \rho_s \left( \sum_{l \in L} \sum_{c \in C} \left\{ \lambda_{lc} \left[ (1 - q_{ls}) p_{lcs} \right] + o_{lc} k_{lcs} + \gamma_{lc} d_{lcs} \right\} + n_l z_{ls} \right)$     (3)

subject to

$$\sum_{c \in C} \left[ (1 - q_{ls}) p_{lcs} + k_{lcs} + d_{lcs} \right] \leq \kappa_l x_l, \quad \forall l \in L, s \in S \tag{4}$$

$$k_{lcs} + d_{lcs} = q_{ls} p_{lcs}, \quad \forall c \in C, l \in L, s \in S \tag{5}$$

$$k_{lcs} - (r_{ls}) p_{lcs} \leq M(1 - z_{ls}), \quad \forall c \in C, l \in L, s \in S \tag{6}$$

$$d_{lcs} - (q_{ls} - r_{ls}) p_{lcs} \leq M(1 - z_{ls}), \quad \forall c \in C, l \in L, s \in S \tag{7}$$

$$d_{lcs} \leq M(z_{ls}), \quad \forall c \in C, l \in L, s \in S \tag{8}$$

$$\sum_{l \in L} \left[ (1 - q_{ls}) p_{lcs} + k_{lcs} \right] = b_c, \quad \forall c \in C, s \in S \tag{9}$$

$$z_{ls} \leq x_l, \quad \forall l \in L, s \in S \tag{10}$$

$$k_{lcs}, d_{lcs}, p_{lcs} \geq 0, \quad \forall c \in C, l \in L, s \in S \tag{11}$$

$$z_{ls} \in \{0, 1\}, \quad \forall l \in L, s \in S \tag{12}$$

$$x_l \in \{0, 1\}, \quad \forall l \in L \tag{13}$$

The objective function (3) consists of five terms. The first term is the fixed cost of selecting facility $l$. The second term represents the expected transportation cost of shipping untainted products. The third and the fourth terms represent the cost of supplying tainted products for the consumers and the cost of discarding tainted products, respectively. Finally, the last term is the cost of inspection which is performing at a facility site.

Constraint set (4) requires a facility to be open if any portion of consumer demand is served from the facility. In addition, it ensures that the total consumer demand assigned to any facility cannot exceed the facility's capacity. Constraint sets (5) through (8) together calculate the amount of tainted product that is shipped to the consumer. Hence, without inspection (when $z_{ls} = 0$), constraint set (8) implies that $d_{lcs} = 0$ and given constraint set (5), all of the tainted goods may reach the consumer. However if inspection is implemented, constraint sets (6) and (7) imply that only products passing inspection will be shipped to the consumer. Constraint set (9) requires that that the demand of every consumer be met. Constraint set (10) implies that inspection can be applied to only the selected set of facilities. Constraint set (11) require that be $k_{lcs}, d_{lcs}$ and $p_{lcs}$ are positive values, and constraint sets (12) and (13) place binary restrictions on variables $z_{ls}$ and $x_l$.

## 3. Generation of Test Data

Let 0 indicate the state that a facility is in an ideal situation and is perfectly working; and 1 indicate the state that a facility is not working with full capacity or the facility works partially. Let's $\Theta_l \in [0.50, 0.95]$ be the probability of facility $l$ being in state 0. Therefore, the assumption that all facilities have an identical probability of working or failing is relaxed [15]. It is assumed in state 0 that there is no tainted product produced at the facility. If a facility is in state 1, the proportion of tainted product is in the range of $[0.10, 0.30]$. The proportion of tainted product that is detected after inspection is in the range of $[0.01, 0.09]$.

To determine the probability of scenario $s$ ($\rho_s$), we need to define a scenario. A scenario is defined as an event where a subset of facilities say $L'$, are in state 0, and facilities in the set $L \setminus L'$ are in state 1. Given the number of facilities $|L|$, the total number of scenarios that facilities can be in state 1 is given by $\sum_{i=1}^{|L|} \binom{|L|}{i} = 2^{|L|} - 1$. If we also include the scenario in which all facilities are in state 0, then the total number of scenarios becomes $2^{|L|}$. Hence, the probability of realizing a scenario $s \in S$ is defined as $\rho_s = \prod_{l \in L'} \Theta_l \prod_{l \in L \setminus L'} (1 - \Theta_l)$. We list other assumptions as follows:

- The fixed cost of opening a facility is drawn from a discrete uniform distribution between $1,000,000 and $2,000,000.
- Demand for each consumer is drawn from a discrete uniform distribution between 100 and 300 units.
- The cost of inspection at each facility is drawn from a discrete uniform distribution between $50,000 and $100,000.
- The cost of shipping untainted products is drawn from a discrete uniform distribution between $100 and $1000.
- The cost of shipping tainted products is drawn from a discrete uniform distribution between $10,000 and $25,000.

- The cost to discard is equal to 25% of the cost of shipping untainted goods.
- The cost of selecting a facility is correlated with the capacity such that the highest capacity has the highest selecting cost.
- The cost of inspection is correlated to the percentage improvement which is the difference between $q_l$ and $r_l$.
- Total capacity is tight and is 30 percent higher than the total demand before performing inspection and discarding tainted items.

## 4. Solution Procedure

### 4.1. Heuristics

We present a few constructive (greedy) heuristics in this section. In constructive algorithms, we start from scratch (empty solution) and construct a solution by assigning values to one decision variable at a time, until a complete solution is generated [16]. Constructive algorithms are popular techniques as they are simple to design. Moreover, their complexity compared to other algorithms such as iterative algorithms is low. However, in most optimization problems, the performance of constructive algorithms may be low as well. Therefore, we also develop improvement algorithms to improve the quality of the solution achieved by constructive algorithms. In our improvement algorithm, we start with a complete solution (i.e., a constructive algorithm solution) and transform it at each iteration using some search operators to hopefully find a better solution.

In our solution procedure, we first determine the set of selected facilities, $x$. Given the set of selected facilities, we determine the values for inspection i.e. $z$. Having $x_l$ and $z_{ls}$ determined and fixed to their binary values, equation (3) reduces to a capacitated transportation problem, which is relatively easier to solve. We call this model SCD-Sub, and its formulation is stated as follows:

$$[\text{SCD-Sub}] \quad \min \sum_{s \in S} \rho_s \left( \sum_{l \in L} \sum_{c \in C} \left\{ \lambda_{lc} \left[ (1 - q_{ls}) p_{lcs} \right] + o_{lc} k_{lcs} + \gamma_{lc} d_{lcs} \right\} \right)$$

Notice that $k_{lcs}$ and $d_{lcs}$ are auxiliary decision variable which depend solely on $p_{lcs}$ and $z_{ls}$. However, given the fact that $z_{ls}$ is already determined and fixed, we can rewrite the SCD-Sub as follows:

$$[\text{SCD-Sub}] \quad \min \sum_{s \in S} \rho_s \left( \sum_{l \in L} \sum_{c \in C} \left\{ \lambda_{lc} (1 - q_{ls}) + o_{lc} q_{ls} (1 - \overline{z}_{ls}) + o_{lc} r_{ls} \overline{z}_{ls} + \gamma_{lc} (1 - \overline{z}_{ls})(q_{ls} - r_{ls}) \right\} p_{lcs} \right)$$
(14)

$$\text{subject to} \quad k_{lcs} + d_{lcs} = q_{ls} p_{lcs}, \quad \forall c \in C, l \in L, s \in S \tag{15}$$

$$k_{lcs} - (r_{ls}) p_{lcs} \leq M(1 - \overline{z}_{ls}), \quad \forall c \in C, l \in L, s \in S \tag{16}$$

$$d_{lcs} - (q_{ls} - r_{ls}) p_{lcs} \leq M(1 - \overline{z}_{ls}), \quad \forall c \in C, l \in L, s \in S \tag{17}$$

$$d_{lcs} \leq M(\overline{z}_{ls}), \quad \forall c \in C, l \in L, s \in S \tag{18}$$

$$\sum_{l \in L} \left[ (1 - q_{ls}) p_{lcs} + k_{lcs} \right] = b_c, \quad \forall c \in C, s \in S \tag{19}$$

$$k_{lcs}, d_{lcs}, p_{lcs} \geq 0, \quad \forall c \in C, l \in L, s \in S \tag{20}$$

where $\overline{z}_{ls}$ is the fixed and known value of $z_{ls}$. We will refer to SCD-Sub in the following.

### 4.1.1. Constructive Heuristics

Our constructive heuristics operate in three stages. In stage one we determine the set of selected facilities $x$. In stage two we determine the set of $z$ (the inspection policy), and finally in the last stage, we solve SCD-Sub in two phases. In the following, all three stages are presented.

Stage 1 Methods

We develop three constructive heuristics to determine the set of selected facilities as follows:

- Basic Greedy Heuristic (BGH): One way to determine vector $x$ is to simply open all the facilities. Therefore, we have $x_l = 1, \forall l \in L$.

- Selective Greedy Heuristic (SGH): In this method, we first start with an empty set for selected facilities. Steps are illustrated in Figure 4.

> 1. Calculate the expected available capacity of a facility. Assume that a facility operates with full capacity with probability $\chi_l$ and let $\bar{\kappa}_l$ denote the random variable for the available capacity where $p(\bar{\kappa}_l = \kappa_l) = \chi_l$ and $p(\bar{\kappa}_l = (1-q_l)\kappa_l) = (1-\chi_l)$. Therefore, the expected available capacity of facility l can be defined as $E(\bar{\kappa}_l) = \chi_l \kappa_l + (1-\chi_l)(1-q_l)\kappa_l$. For the sake of simplicity, we consider $\chi_l = \frac{1}{2}$ in our computations.
> 2. Evaluate the total cost of selecting a facility. The estimated total costs of selecting a facility includes the fixed cost, and mean costs of shipping untainted products, shipping tainted products, and discarding tainted products. These costs are computed as $\frac{\sum_{c=1}^{|C|} \lambda_{lc}}{|C|}$, $\frac{\sum_{c=1}^{|C|} o_{lc}}{|C|}$, and $\frac{\sum_{c=1}^{|C|} \gamma_{lc}}{|C|}$, respectively.
> 3. Sort the calculated values in an increasing order, obtain the sequence of the facilities, and locate them in the set $E$. Let $\kappa_{[i]}$ be the corresponding expected available capacity of $i^{th}$ facility in set $E$.
> 4. $i \leftarrow 1$.
> 5. Choose the $i^{th}$ cheapest facility from the set $E$.
> 6. **if** $\sum_{j=1}^{i} \kappa_{[j]} \geq \sum_{c \in C} b_c$, go to 7, **else** $i \leftarrow i+1$ and go to 5.
> 7. **end**

Figure 4 Pseudocode of the Selective Greedy Heuristic (SGH)

- Capacity-Based Greedy Heuristic (CBGH): In this method, we first start with an empty set for selected facilities. Then we choose a facility from the set of remaining candidates that reduces the total demand of the consumer ($\sum_{c \in C} b_c$) the most. Steps are illustrated in Figure 5.

> 1. Calculate the expected available capacity as presented in step 1 of SGH.
> 2. Sort the calculated values in step 1 in a decreasing order, obtain the sequence of the facilities, and locate them in the set H. Let $\kappa_{[i]}$ be the corresponding expected available capacity of $i^{th}$ facility in set H.
> 3. Choose the $i^{th}$ facility from the set H until satisfying $\sum_{j=1}^{i} \kappa_{[j]} \geq \sum_{c \in C} b_c$.

Figure 5 Pseudocode of the Capacity-Based Greedy Heuristic (CBGH):

Stage 2 Methods

Once the selected facilities have been set, we determine the inspection policy ($z_{ls}$) for each scenario. We consider the following construction heuristics:

- Failed Scenario Inspection Heuristic (FSIH): Consider $S'_l, \forall l \in L$ as the set of scenarios where facility $l$ works with full capacity and $S''_l, \forall l \in L$ as the set of scenarios where facility $l$ will produce tainted products where $S'_l \cup S''_l = S$. We define $z_{ls}$ as follows: $z_{ls} = \begin{cases} 1, & s \in S_l'' \wedge x_l = 1, \forall l \in L, \forall s \in S \\ 0, & s \in S_l' \vee x_l = 0, \forall l \in L, \forall s \in S \end{cases}$. This heuristic performs inspection for only those facilities that are selected in stage 1 and belong to set of scenarios where facility $l$ produces tainted products.

- Greedy Inspection Heuristic (GIH): In this method, we define a desired shipping untainted level $\Delta \left(\Delta \in (0,1]\right)$, where $(\Delta)100\%$ of the shipping products to consumers must be untainted. Let's start with an empty set for $z$. Given a scenario $s$ where $s \in S_l''$, following relation should be satisfied $\sum_{l \in L}\left[q_{ls}x_l\kappa_l(1-z_{ls}) + r_{ls}x_l\kappa_l z_{ls}\right] \le (1-\Delta)\sum_{c \in C}b_c$. Otherwise, we perform inspection until we reach the desired level of untainted products. Note that, we consider facilities in decreasing order of maximum reduction in the fraction tainted products or $\max\{(q_l - r_l) \mid \forall l \in L\}$. We consider $\Delta = 0.90$ in our computations.

- Random Greedy Inspection Heuristic (RGIH): The basic idea of this method is to estimate how much we can save by implementing inspection in a facility. If this saving is significant enough, then the inspection for a facility will be implemented. Given $x$ as the set of selected facilities, the steps are defined in Figure 6.

**for** s = 1 to |S|
1. Determine the amount of saving: calculate the saving for each facility by following equation. $\phi_{ls} = \dfrac{(q_{ls} - r_{ls})\sum_{c=1}^{|C|} o_{lc}}{n_l + (q_{ls} - r_{ls})\sum_{c=1}^{|C|} \gamma_{lc}}$, $\forall l \in L$, where $\phi_{ls}$ is the amount of saving for facility $l$ in scenario $s$. Let the vector $\Phi$ be the set of calculated savings where $|\Phi| = |x|$.
2. Normalize vector $\Phi$: $\hat{\Phi} = \dfrac{\Phi}{\|\Phi\|}$.
3. Generate a set of random numbers, $\hat{R} \in (0,1)$, where $|\hat{R}| = |\Phi| = |x|$.
4. Compare each element of vector $\hat{R}$ with the corresponding element of vector $\hat{\Phi}$. If that is greater than the corresponding normalized saving value, then $z_{ls} = 0$, and 1 otherwise.
**end for**

Figure 6 Pseudocode of the Random Greedy Inspection Heuristic (RGIH)

## Stage 3 Methods

In this stage we solve SCD-Sub in two phases. In the first phase, we set $p_{lcs}$ in a greedy fashion, based on the unit transportation cost to the consumers ($\lambda_{lc}$) and capacity of the selected facilities. Note that in Figure 7, $a_c$ represents the demand of consumer $c$, and $g_l$ represents the capacity for facility $l$. In the second phase, given the obtained value for $p_{lcs}$, we simply compute the values of auxiliary variables $k_{lcs}$ and $d_{lcs}$.

---

**Input:** $c \in C$ set of consumers, $l \in L$ set of facilities, $s \in S$ set of scenarios, $f_l$, $\kappa_l, n_l$, $b_c$, $\rho_s$, $\lambda_{lc}$, $o_{lc}$, $q_{ls}$ and $r_{ls}$.

**Output:** $x_l, z_{ls}, p_{lcs}, k_{lcs}, d_{lcs}$ and the total cost.

**Stage One:** Determine $x$

1. Use *BGH, SGH* or *CBGH*.

**Stage Two:** Determine $z_{ls}$

2. Use *FSIH, GIH,* or *RGIH*.

**Stage Three:** Solve SCD-Sub

**Phase 1**

3.     $\forall s \in S$
4.     Sort $\lambda_{lc}$ in increasing order. $a_c \leftarrow b_c$ ; $g_l \leftarrow \kappa_l$.
5.        $\forall c \in C$
6.           **while** $a_c > 0$ **do**
7.              $l = \arg\min\{\lambda_{lc}\}\ \forall l \in L, c \in C$
8.              **if** $g_l > a_c$ **then**
9.                  $p_{lcs} \leftarrow a_c, g_l \leftarrow g_l - a_c, a_c \leftarrow 0$
10.             **else if** $g_l > 0$ **then**
11.                 $p_{lcs} \leftarrow g_l, a_c \leftarrow a_c - g_l, g_l \leftarrow 0$
12.             **end if**
13.          **end while**
14.        **end for**
15. **end for**

**Phase 2**

16. $\forall s \in S, \forall l \in L$
17.     **if** $z_{ls} = 0$ **then**
18.        $\forall c \in C$
19.           $d_{lcs} \leftarrow 0$
20.           $k_{lcs} \leftarrow q_{ls} p_{lcs}$
21.        **end for**
22.     **else if** $Z_{ls} = 1$ **then**
23.        $\forall c \in C$
24.           $d_{lcs} \leftarrow (q_{ls} - r_{ls}) p_{lcs}$
25.           $k_{lcs} \leftarrow r_{ls} p_{lcs}$
26.        **end for**
27.     **end if**
28. **end for**
29. **end**

Figure 7 Pseudocode to Solve Problem SCD-Sub

## 4.2. Improvement Heuristics

In this section we develop improvement heuristics to improve the solution obtained from one of the heuristic methods presented above (note that improvement heuristics operate on a solution found by a constructive heuristic.) First, we present an improvement heuristic that begins with a feasible solution and seeks to improve upon it. The improvement heuristic iteratively closes one facility if the facility is already selected and opens a facility if a facility is not selected. This iteration enables us to generate a new neighborhood and explore if the new set of selected facilities provides a cheaper solution or not. In order to maintain feasibility, only moves are allowed which provide enough capacity to satisfy the total demand of the consumers. The details of this heuristic are shown in Figure 8.

1. **for** each facility $l \in x$, if $x_l = 1$ then $x_l \leftarrow 0$ otherwise **if** $x_l = 0$ **then** $x_l \leftarrow 1$. Let $\psi$ be the new set of selected facilities ($|\psi| = |x|$).
2. Compute saving as:
$$\sigma_l = SCD(x) - SCD(\psi)$$
3. **If** $\sigma_l < 0$ **then** $x_l \leftarrow \psi_l$. Go to 1.
4. **end for**

Figure 8 Pseudocode of the local_$x$

In the second improvement algorithm, we apply a Variable Neighborhood Search (VNS). The basic idea of VNS is to find a set of predefined neighborhoods to achieve a better solution. It explores either at random or deterministically a set of neighborhoods to get different local optima and to escape from local optima (for general pseudocode of VNS see [17]). The purpose of the second improvement is to minimize transportation cost for each individual scenario, i.e., minimize $\sum_{l \in L}\sum_{c \in C}\left(\lambda_{lc}\left[(1-q_{ls})p_{lcs}\right]\right)$ (21). We implement the VNS for our problem as follows:

1. Set $k \leftarrow 1$
2. **for** all $s \in S$
3.     while $k < K_{max}$
**Shaking**:
4.     Set $u \leftarrow p$
5.     Define a neighborhood strategy in u
6.     Apply a mechanism to generate a new solution for u
**Improve or not:**
7.     Calculate the cost for $u$ from equation (21).
8.     If $cost(u) < cost(p)$ then $p \leftarrow u$ else $k \leftarrow k+1$
9.    **end while**
10. **end for**

Figure 9 VNS for improving the transportation cost

The neighborhood strategy that we apply is structured by randomly choosing a point in the matrix of transportation. Subsequently, we identify the closed path leading to that point which consists of horizontal and vertical lines as illustrated. In order to generate a new solution, we move $\hat{R}$ unit(s) from the chosen point and another point at a corner of the closed path and modify the remaining points at the other corners of the closed path to reflect this move. Note that $\hat{R}$ is a random variable over the range of zero and the minimum value of the four selected points. This scheme is demonstrated in Figure 10. The selected point is shown by ∗.

| Consumer / facility | 1 | 2 | 3 | 4 |
|---|---|---|---|---|
| 1 | 50 | 60 + | 85 | 85 - |
| 2 | 21 | 45 | 29 | 29 |
| 3 | 81 | 36* - | 73 | 73 + |
| 4 | 62 | 78 | 91 | 20 |

Figure 10 Neighborhood strategy in VNS

### 4.3. Simulated Annealing

Simulated Annealing (SA) is a metaheuristic approach inspired by nature. In this case, the process of a heated metal being cooled at a controlled rate (annealed) to improve its physical properties is simulated. The method was popularized by the work of Kirkpatrick et al. [18] which continued the earlier work of Metropolis et al. [19]. The fundamental idea is to allow moves resulting in solutions of worse quality than the current solution in order to escape from local optima [16]. The probability of doing such a move is decreased during the search.

An important consideration in SA is to set the initial value of the temperature $T_0$. If the initial temperature is set very high, the search may be relatively close to a random local search. Otherwise, if the initial temperature is very low, the search might degenerate to an improving local search algorithm [16] as the probability of accepting worse moves decreases too quickly. Another important factor is cooling schedule. The choice of a suitable cooling schedule is crucial for the performance of the algorithm. The cooling schedule defines the value of temperature $T$ at every iteration. Figure 11 outlines the implementation of SA in more detail.

### 4.3.1. Defining Initial Temperature and Cooling Schedule

The temperature $T$ is decreased during the search process, thus at the beginning of the search the probability of accepting uphill moves is high and it gradually decreases. As stated, the choice of an appropriate cooling schedule and initial value of temperature are crucial for the performance of the algorithm. The cooling schedule defines the value of $T$ at each iteration $k$, $T_{k+1} = R(T_k, k)$, where $R(T_k, k)$ is a function of the temperature at the previous step and of the iteration number. In this paper, we use one of the most common cooling schedule which follows a geometric law as $T_{k+1} = \theta T_k$, where $\theta \in (0,1)$, which corresponds to an exponential decay of the temperature [20]. Furthermore, experience has shown that $\theta$ should be between 0.5 and 0.99 (see [16]). Hence, we considered four values, $\theta = 0.95, 0.90, 0.80$ and $0.75$; and we obtained the best minimum regret in less computational time at $\theta = 0.75$.

Another important factor in SA is to define the initial value of the temperature $T_0$ properly. There is a tradeoff between a very high initial temperature and a lower one. The high temperature explores more of the solution space at the cost of increased running time. For this research, we use acceptance deviation methods proposed by Huang et al. [21]. The starting temperature is computed by $t\sigma$ using preliminary experimentations on each data set, where σ represents the standard deviation of difference between values of objective functions and $t = {-3}/{ln(\varpi)}$ with the acceptance probability of $\varpi$, which is greater than 3σ. Finally, a sufficient number of iterations at each temperature should be performed. If too few iterations are performed at each temperature, the algorithm may not be able to reach the global optimum. Given the presented formula and after several experiments, we set the value for the initial temperature, $T_0 = 8000$.

### 4.3.2. Neighborhood Selection

The manner in which a metaheuristic technique moves from one solution to its neighbor is a critical component. In our SA algorithm, we define a neighborhood which combines four neighborhood structures: (1) swapping one randomly selected facility with another randomly selected facility (SA-$swap$), (2) selecting one more facility (SA-$add$), (3) closing one selected facility (SA-$remove$), and finally (4) closing two facilities while selecting another two (SA-2$swap$). Note that we apply the same neighbor strategy to determine the $z$ and afterward compute the values of $p_{lcs}^{iter}$, $k_{lcs}^{iter}$, and $d_{lcs}^{iter}$ by using Constructive 1.

### 4.3.3. Stopping Criterion

Various stopping criteria have been developed in the literature. A popular stopping criteria and the temperature reaches a set value (such as 0.01). Another criterion can be completing a predetermined number of iterations. In this paper, a combination of these two criteria is considered in which we stop at the earlier of the temperature reaching 0.01 or the completion of 100 (350) iterations for small (large) size problems.

1. Initialize the parameters of the annealing schedule (Initial temperature, final temperature and total number of iterations)
2. Generate an initial solution by determining vector $x^0$, $z^0$, $p_{lcs}^0$, $k_{lcs}^0$, $d_{lcs}^0$ by the represented constructive or improvement heuristics and define relevant total cost $f(x^0, z^0, p_{lcs}^0, k_{lcs}^0, d_{lcs}^0)$
3. $iter \leftarrow 1$; Temperature $\leftarrow$ Initial Temperature
4. **while** Temperature > Final Temperature or $iter$ < total number of iterations **do**
5.     **while** done=false
6.         aZeroElem $\leftarrow$ Number of zero elements in vector $x^{iter}$ and $z^{iter}$
7.         aOneElem $\leftarrow$ Number of one elements in vector $x^{iter}$ and $z^{iter}$
8.         aRand $\leftarrow$ Generate a Random Number
9.         **if** $0 \le aRand < \frac{1}{4}$
10.             create a new solution using SA-*swap* method and return $x^{iter}$ and $z^{iter}$
11.             done $\leftarrow$ true
12.         **else if** $\frac{1}{4} \le aRand < \frac{1}{2}$ **and** $aZeroElem \ge 1$
13.             create a new solution using SA-*add* method and return $x^{iter}$ and $z^{iter}$
14.             done $\leftarrow$ true
15.         **else if** $\frac{1}{2} \le aRand < \frac{3}{4}$ **and** $aOneElem > 1$
16.             create a new solution using SA-*remove* method and return $x^{iter}$ and $z^{iter}$
17.             done $\leftarrow$ true
18.         **else if** $\frac{3}{4} \le aRand < 1$ **and** $aOneElem \ge 2$
19.             create a new solution using SA-2*swap* method and return $x^{iter}$ and $z^{iter}$
20.             done $\leftarrow$ true
21.     **end while**
22. Obtain the values of $p_{lcs}^{iter}$, $k_{lcs}^{iter}$, and $d_{lcs}^{iter}$ by using SCD-Sub.
23. **if** $f(x^{iter}) - f(x^0) \le 0$ **then** $f(x^0) = f(x^{iter})$, $x^0 = x^{iter}$
24. Update Temperature
25. **end while**
26. return the final solution

Figure 11 Pseudo code of the SA algorithm

### 4.4. Commercial Software

The optimization problem is modeled with the AMPL mathematical programming language and solved with Gurobi version 4.5.6. Each problem instance is solved on 4 cores (threads=4) of a Dell Optiplex 980 with an Intel Core i7 860 Quad @ 2.80GHz, and 16GB RAM. The operating system is Windows 7 Enterprise 64 bit. In our computational analysis, we terminate Gurobi when the CPU time limit of 14400 seconds is reached. Table 1 summarizes the results from the solution, and the discussion is presented in section 5.

### 5. Computational Results

In this section, we perform computational experiments to assess the effectiveness of the algorithms. In section 4.1.1, we presented three heuristics (BGH, SGH, and CBGH) to determine set of selected facilities, $x$, and also three heuristics (FSIH, GIH, and RGIH) to determine the set of inspections to conduct, $z$. By combining these six heuristics, we construct nine different heuristic for determining $x$ and $z$. For instance, our first heuristic can be denoted as BGH&FSIH. Finally, we employ the greedy heuristic presented in Figure 7 to solve SCD-Sub. All the algorithms were created and executed in MATLAB 7.9 (2009b) and tested on a single core of a Dell OptiPlex 980 computer running the Windows 7 Enterprise 64 bit operating system with an Intel(R) Core(TM) i7 CPU860@ 2.80GHz, and 8GB RAM.

We consider 12 sets of problems with 10 data instances in each. Hence, we solve in total 120 instances of varying sizes as illustrated in Table 1. The second, third and the fourth columns represent the size of the problems under consideration. We also report the average of the optimal value and average solution time for each set. Finally, the last column represents the total number of optimal solutions obtained from 10 data instances.

Table 1 Test problems' sizes and the corresponding optimum solution values and times

| Set no. | No. of consumers | No. of facilities | No. of scenarios | Avg. optimal value/best values found | Avg. optimum time (s) | No. of optimal solutions in 10 instances |
|---|---|---|---|---|---|---|
| 1 | 2 | 2 | 4 | 3132033.5 | 0.014 | 10 |
| 2 | 2 | 5 | 32 | 6031316.5 | 17.6 | 10 |
| 3 | 2 | 10 | 1024 | 12190855.2* | 11267.6 | 5 |
| 4 | 5 | 2 | 4 | 3857839.1 | 0.02 | 10 |
| 5 | 5 | 5 | 32 | 7130036.5 | 62.9 | 10 |
| 6 | 5 | 10 | 1024 | 12632921.9* | 13102.7 | 1 |
| 7 | 10 | 2 | 4 | 4498383.3 | 0.022 | 10 |
| 8 | 10 | 5 | 32 | 7248655.8 | 1504.8 | 9 |
| 9 | 10 | 10 | 1024 | 13207650.1* | ** | 0 |
| 10 | 20 | 2 | 4 | 5526660.8 | 0.03 | 10 |
| 11 | 20 | 5 | 32 | 8503834.6 | 325.7 | 10 |
| 12 | 20 | 10 | 1024 | 15095012.3* | ** | 0 |

*: Average of best objective values found

**: The CPU time exceeded the prescribed time limit of 14400 seconds.

As observed from Table 1, increasing the number of facilities implies an increase in the number of scenarios and the size of the problem has a sigficant impact on the solution time. For instance, for the case of 10 consumers or 20 consumers and 10 facilities, Gurobi did not return any optimal solutions within the prescribed time limit of 14400 seconds. In order to calculate the relative optimality gap, we use the objective function value that is provided by Gurobi when the prescribed time limit is reached. To assess each heuristic, we consider solution quality and solution (computational) time. For the solution quality, we consider a quality criterion which is the gap between the result of heuristic/SA and the optimal/best solution obtained from Gurobi. This gap is defined according to the following equation:

$$\% \text{gap} = \frac{\left(SA(or\ Heuristics)\ Solution - Best\ Found\ (or\ Optimal\ Solution)\right)}{Best\ Found\ (or\ Optimal\ Solution)}$$

Furthermore, given the randomness characteristic of GIH, RGIH, Local_$x$, VNS_$P$, and SA, the corresponding objective values and solution times are the average across thirty independent replications. Table 2 reports the result for 2 facilities and 2, 5, 10, and 20 consumers. Note that bold-faced values indicate achievement of the best optimality among constructive heuristics and improvement heuristics/SA, respectively.

The results in Table 2 show that, regardless of the number of consumers, the heuristic algorithms always provide solutions within 3% of the solution found by Gurobi. Heuristic algorithms are fast and their solution time generally does not vary with the number of consumers. We see that SGH&FSIH and CBGH&FSIH algorithms provide better solution quality and Local_$x$ does not provide any improvement in the solution of the constructive heuristics. The VNS_$P$ procedure provides better quality solutions than the Local_$x$, however, this improvement comes with an increase in the solution time. Another observation from Table 2 is that even though the solution time for SA algorithm is notably higher than the other algorithms, its performance is not as good as VNS_$P$ when we limit the problem instances to 2 facilities.

Table 2 Comparison of algorithms results for 2 facilities

| | | Constructive Heuristics | | | | | | | | | Improv. Heuristic and Metaheuristic | | |
|---|---|---|---|---|---|---|---|---|---|---|---|---|---|
| | | BGH&FSIH | BGH&GIH | BGH&RGIH | SGH&FSIH | SGH&GIH | SGH&RGIH | CBGH&FSIH | CBGH&GIH | CBGH&RGIH | Local_$x$ | VNS_$P$ | SA |
| 2 Consumers | Min gap | 0.00% | 0.00% | 0.00% | 0.00% | 0.00% | 0.00% | 0.00% | 0.08% | 0.00% | 0.00% | 0.00% | 0.00% |
| | Avg. gap | 0.40% | 0.47% | 0.22% | **0.21%** | 0.64% | 0.36% | 0.36% | 0.56% | 0.25% | 0.21% | **0.06%** | 0.25% |
| | Max gap | 0.89% | 0.77% | 1.59% | 0.76% | 2.49% | 0.76% | 0.56% | 1.59% | 0.59% | 0.56% | 0.09% | 0.89% |
| | Avg. time (s) | 0.002 | 0.002 | 0.002 | 0.002 | 0.002 | 0.002 | 0.002 | 0.002 | 0.002 | 0.008 | 0.021 | 0.047 |
| 5 Consumers | Min gap | 0.00% | 0.00% | 0.00% | 0.00% | 0.00% | 0.00% | 0.00% | 0.00% | 0.00% | 0.00% | 0.00% | 0.00% |
| | Avg. gap | 1.34% | 1.59% | 1.86% | 0.97% | 1.36% | 1.33% | **0.84%** | 1.08% | 1.23% | 0.84% | **0.25%** | 0.55% |
| | Max gap | 3.75% | 3.62% | 3.83% | 3.83% | 3.83% | 3.83% | 1.47% | 2.53% | 2.33% | 1.47% | 1.37% | 1.47% |
| | Avg. time (s) | 0.002 | 0.002 | 0.002 | 0.002 | 0.002 | 0.002 | 0.002 | 0.002 | 0.002 | 0.006 | 0.044 | 0.047 |
| 10 Consumers | Min gap | 0.00% | 0.00% | 0.00% | 0.00% | 0.00% | 0.00% | 0.00% | 0.00% | 0.00% | 0.00% | 0.00% | 0.00% |
| | Avg. gap | 0.96% | 2.12% | 0.68% | **0.48%** | 2.19% | 1.08% | **0.48%** | 0.81% | 1.08% | 0.48% | **0.05%** | 0.50% |
| | Max gap | 4.77% | 7.02% | 1.84% | 1.84% | 7.02% | 4.55% | 1.84% | 2.16% | 2.84% | 1.84% | 0.40% | 1.84% |
| | Avg. time (s) | 0.003 | 0.003 | 0.003 | 0.003 | 0.003 | 0.003 | 0.003 | 0.003 | 0.003 | 0.006 | 0.03 | 0.14 |
| 20 Consumers | Min gap | 0.00% | 0.00% | 0.43% | 0.00% | 0.00% | 0.00% | 0.00% | 0.00% | 0.00% | 0.00% | 0.00% | 0.00% |
| | Avg. gap | **1.58%** | **1.58%** | 2.31% | **1.58%** | **1.58%** | 1.97% | **1.58%** | 1.90% | 1.90% | 1.58% | **1.29%** | 1.58% |
| | Max gap | 5.37% | 5.37% | 5.37% | 5.37% | 5.37% | 5.37% | 5.37% | 5.37% | 5.37% | 5.37% | 4.49% | 5.37% |
| | Avg. time (s) | 0.002 | 0.002 | 0.002 | 0.002 | 0.002 | 0.002 | 0.002 | 0.002 | 0.002 | 0.005 | 0.023 | 0.117 |

Bold-faced values indicate achievement of the best optimality.

We now evaluate the effectiveness of our algorithms for five facilities. The results are presented in Table 3. The performance of BGH&FSIH, BGH&GIH, and BGH&RGIH is not good. The reason is that in these three heuristics we use BGH to select all the facilities while the total demand can be satisfied by selecting fewer facilities. SGH&FSIH, SGH&GIH, and SGH&RGIH provide solutions on average within 8% of the best found solution with a remarkably fast solution time in comparison to the optimal solution time. Both Local_$x$ and VNS_$P$ are capable of improving the solution quality even for a larger number of consumers and the average solution gap is within 5% of the optimal solution. In particular, SA clearly provides the best overall solution cost for the range of problems tested and requires only a moderate extra computational time than other algorithms. SA achieves solutions which are in average within 3% of the optimality gap. For 5 facilities and $|C| \in \{10\}$, as shown in Table 1, we found 9 optimal solutions in 10 data instances. Therefore, in Table 3 we show the gap with the optimal and non-optimal solutions (or best solutions found) individually.

Table 3 Comparison of algorithms results for 5 facilities

| | | Constructive Heuristics | | | | | | | | | Improv. Heuristic and Metaheuristic | | |
|---|---|---|---|---|---|---|---|---|---|---|---|---|---|
| | | BGH&FSIH | BGH&GIH | BGH&RGIH | SGH&FSIH | SGH&GIH | SGH&RGIH | CBGH&FSIH | CBGH&GIH | CBGH&RGIH | Local_x | VNS_P | SA |
| 2 Consumers | Min gap* | 7.09% | 7.09% | 7.09% | 0.52% | 1.34% | 1.34% | 0.29% | 0.29% | 0.30% | 0.29% | 0.26% | 0.14% |
| | Avg. gap* | 16.05% | 15.68% | 15.74% | 3.83% | 3.99% | 3.98% | **2.37%** | 2.61% | 2.49% | 2.16% | 1.92% | **0.89%** |
| | Max gap* | 30.08% | 30.08% | 28.65% | 12.36% | 12.37% | 12.32% | 7.46% | 7.53% | 7.70% | 6.40% | 2.02% | 1.08% |
| | Avg. time (s) | 0.002 | 0.002 | 0.002 | 0.002 | 0.002 | 0.002 | 0.002 | 0.002 | 0.002 | 0.01 | 0.066 | 0.146 |
| 5 Consumers | Min gap* | 13.10% | 14.38% | 13.74% | 1.00% | 0.98% | 1.45% | 1.44% | 2.38% | 1.47% | 0.58% | 0.00% | 0.04% |
| | Avg. gap* | 18.37% | 18.35% | 18.30% | **4.13%** | 4.20% | 4.11% | 4.52% | 4.48% | 4.38% | 1.99% | 1.23% | **1.18%** |
| | Max gap* | 25.37% | 25.33% | 25.64% | 7.78% | 7.60% | 7.54% | 11.14% | 11.11% | 11.27% | 3.50% | 2.83% | 2.66% |
| | Avg. time (s) | 0.002 | 0.002 | 0.002 | 0.002 | 0.002 | 0.003 | 0.002 | 0.002 | 0.002 | 0.009 | 0.173 | 0.208 |
| 10 Consumers | Min gap* | 8.37% | 9.70% | 9.09% | 1.50% | 1.28% | 3.24% | 3.59% | 3.69% | 3.48% | 1.44% | 0.98% | 0.43% |
| | Avg. gap* | 14.24% | 14.51% | 14.90% | **4.02%** | 4.41% | 4.72% | 7.81% | 7.72% | 7.37% | 2.71% | 2.21% | **1.49%** |
| | Max gap* | 19.20% | 18.13% | 19.08% | 6.08% | 6.09% | 6.60% | 12.09% | 11.44% | 10.25% | 3.59% | 3.58% | 2.21% |
| | Avg. gap w/ non-opt sol. ** | 9.63% | 9.70% | 9.83% | **2.17%** | 3.74% | 3..86% | 8.31% | 8.70% | 8.15% | 1.03% | 0.03% | *-0.11%* |
| | Avg. time (s) | 0.002 | 0.002 | 0.002 | 0.002 | 0.002 | 0.003 | 0.002 | 0.002 | 0.002 | 0.014 | 0.318 | 0.709 |
| 20 Consumers | Min gap* | 2.01% | 2.01% | 3.30% | 3.12% | 3.11% | 4.02% | 6.86% | 6.87% | 6.94% | 1.97% | 1.89% | 1.75% |
| | Avg. gap* | 9.70% | 9.70% | 10.91% | 6.04% | **6.02%** | 7.94% | 12.75% | 12.83% | 13.78% | 4.65% | 3.53% | **3.17%** |
| | Max gap* | 15.06% | 15.06% | 18.32% | 12.69% | 13.32% | 13.91% | 15.77% | 15.77% | 17.68% | 9.41% | 5.67% | 5.37% |
| | Avg. time (s) | 0.002 | 0.002 | 0.002 | 0.002 | 0.002 | 0.002 | 0.002 | 0.002 | 0.002 | 0.029 | 0.635 | 1.049 |

Bold-faced values indicate achievement of the best optimality/solution gap.

Italicized indicate a better solution than the best solution found by Gurobi within the time limit.

*: Values indicate the average gap with optimal solutions found

**: Values indicate the average gap with the best solution found

Realistic sized problems are commonly larger than those tested above. Hence, we examine a larger size problem for 10 facilities. In Table 1 we show that for set 9 and set 12 we were not able to find the optimal solution for any of the 10 instances in 14400 seconds. In addition, in set 3 and 5 only 50% and 10% of the data instances were solved to optimality, respectively. This indicates how increasing the number of facilities and correspondingly the number of scenarios has a significant impact on the solution time. We presented the result of algorithms for the tested problem in Table 4. Negative values in Table 4 indicate that the heuristics or SA achieved a better solution than the best solution found by Gurobi. For 10 facilities and $|C| \in \{10,20\}$, SGH&FSIH, SGH&GIH, and SGH&RGIH perform well based on the average solution gap. For the case of 10 facilities and $|C| \in \{2,5\}$ consumers CBGH&FSIH, CBGH&GIH, and CBGH&RGIH achieved a better performance. However, the SA solutions outperform those found by all the other algorithms, even though they require less computational time than SA. Also, the minimum and maximum gap is usually somewhat better for the SA. Hence, for large size problems we recommend using the SA algorithm, although reasonable results can still be achieved by some of the algorithms. For 10 facilities and $|C| \in \{5\}$, we found only 1 optimal solution in 10 data instances. Hence, we separate the result for this data instance from the others and display the gap between the optimal solution and the algorithms in the corresponding row of Table 4.

Table 4 Comparison of algorithms results for 10 facilities

| | | Constructive Heuristics | | | | | | | | | Improv. Heuristic and Metaheuristic | | |
|---|---|---|---|---|---|---|---|---|---|---|---|---|---|
| | | BGH&FSIH | BGH&GIH | BGH&RGIH | SGH&FSIH | SGH&GIH | SGH&RGIH | CBGH&FSIH | CBGH&GIH | CBGH&RGIH | Local_x | VNS_P | SA |
| 2 Consumers | Min gap* | 20.01% | 21.25% | 20.52% | 4.24% | 4.17% | 2.52% | 0.46% | 0.46% | 0.43% | 0.13% | 0.09% | 0.09% |
| | Avg. gap* | 24.10% | 24.02% | 22.75% | 7.04% | 6.86% | 6.57% | 3.93% | 3.89% | **3.56%** | 1.14% | 0.54% | **0.12%** |
| | Max gap* | 29.06% | 28.11% | 28.00% | 10.63% | 8.72% | 8.78% | 9.65% | 9.98% | 9.29% | 2.33% | 1.24% | 0.35% |
| | Avg. gap w/ non-opt sol. ** | 21.17% | 21.26% | 21.20% | 6.72% | 6.69% | 6.44% | 3.21% | 3.45% | **3.14%** | 1.09% | 0.47% | ***-0.02%*** |
| | Avg. time (s) | 0.005 | 0.005 | 0.005 | 0.005 | 0.005 | 0.006 | 0.005 | 0.005 | 0.005 | 0.218 | 3.838 | 4.075 |
| 5 Consumers | Avg. gap* | 20.25% | 21.42% | 22.63% | 8.41% | 7.14% | 5.94% | **0.89%** | **0.89%** | 0.90% | 0.89% | 0.87% | **0.85%** |
| | Min gap w/ non-opt. sol. ** | 10.37% | 10.31% | 0.70% | 4.92% | 4.32% | 3.72% | 0.89% | 0.89% | 0.90% | 0.87% | *-0.31%* | *-0.36%* |
| | Avg. gap w/ non-opt. sol. ** | 15.63% | 15.69% | 16.96% | 5.18% | 5.67% | 5.62% | 4.63% | 4.62% | **4.60%** | 1.93% | 0.40% | **0.38%** |
| | Max gap w/ non-opt. sol. ** | 25.23% | 24.98% | 24.70% | 15.40% | 9.52% | 8.89% | 8.67% | 8.83% | 9.78% | 4.22% | 1.40% | 1.73% |
| | Avg. time (s) | 0.007 | 0.007 | 0.007 | 0.007 | 0.008 | 0.008 | 0.007 | 0.007 | 0.008 | 0.196 | 9.068 | 16.907 |
| 10 Consumers | Min gap w/ non-opt. sol. ** | 11.14% | 11.08% | 11.30% | -3.89% | -3.64% | -2.58% | -4.34% | -4.43% | -4.29% | -4.49% | -6.28% | -6.52% |
| | Avg. gap w/ non-opt. sol. ** | 18.25% | 18.36% | 17.75% | **3.38%** | **3.38%** | 4.05% | 5.20% | 5.21% | 5.59% | 0.95% | *-0.37%* | ***-0.88%*** |
| | Max gap w/ non-opt. sol. ** | 25.68% | 25.58% | 23.73% | 8.54% | 9.54% | 10.95% | 14.42% | 12.95% | 12.98% | 5.87% | 4.02% | 1.73% |
| | Avg. time (s) | 0.011 | 0.011 | 0.011 | 0.011 | 0.011 | 0.011 | 0.011 | 0.011 | 0.011 | 0.825 | 15.839 | 19.971 |
| 20 Consumers | Min gap w/ non-opt. sol. ** | -0.89% | -0.83% | -0.15% | -12.02% | -11.99% | -10.70% | -7.63% | -7.59% | -5.52% | -13.33% | -13.33% | -14.28% |
| | Avg. gap w/ non-opt. sol. ** | 11.26% | 11.30% | 11.77% | ***-1.85%*** | *-1.74%* | 0.58% | 4.59% | 4.62% | 5.37% | *-2.34%* | *-2.76%* | ***-3.77%*** |
| | Max gap w/ non-opt. sol. ** | 17.78% | 17.80% | 18.06% | 2.70% | 2.70% | 5.45% | 17.03% | 17.05% | 17.03% | 2.42% | 2.37% | 0.59% |
| | Avg. time (s) | 0.018 | 0.018 | 0.018 | 0.018 | 0.018 | 0.018 | 0.018 | 0.018 | 0.018 | 2.266 | 35.945 | 44.989 |

Bold-faced values indicate achievement of the best optimality/solution gap.
Italicized indicate a better solution than the best solution found by Gurobi within the time limit.

*: Values indicate the average gap with optimal solutions found

**: Values indicate the average gap with the best solution found

Insert Figure 12 somewhere here

Insert Figure 13 somewhere here

It is observable from the result that the number of facilities and consequently the number of scenarios has a significant impact on the computational time in our model. However, the results indicate the effectiveness of the SA algorithm we proposed, particularly for larger sized problems. For problems in practice (that can have even larger sizes), our SA heuristic shows promising results.

## 6. Conclusions and Future Research

In response to some catastrophic events, particularly in healthcare/pharmaceutical supply chains, this research addresses a supply chain network design to hedge against the risk of supply disruptions and sending tainted materials to consumers. We considered a mixed-integer stochastic programming model with capacitated facilities. The model was formulated as a two-stage optimization problem. The aim of the model consists of the facility selection, actual capacity allocation among the consumers, and determination of inspection policy with the objective of minimizing the total cost. The impact of supply/capacity uncertainty is explicitly modeled in all our models in order to design a reliable supply chain network. To capture the uncertainty, a scenario-based approach was presented.

Experience from solving the problem using commercial software indicated that the number of facilities, and consequently the number of scenarios, has a significant impact on the computational time. As a result, we developed several heuristic methods and a metaheuristic approach to effectively solve the presented model.

Based on our computational studies, the SA approach is not efficient in terms of solution quality and solution time for the small size problems or small number of scenarios. However, some of the heuristics, in particular SGH&FSIH, SGH&GIH, SGH&RGIH and CBGH&FSIH, achieved good solution qualities in a more reasonable time when compared to the optimal or best found solution. Local_$x$ and VNS_$P$ were able to improve the solutions obtained from constructive heuristics. Therefore, constructive and improvement heuristics are preferable on small sized problems. However, for practical sized problems, i.e. 10 facilities and more, SA outperforms constructive and improvement heuristics, even though it requires higher computational time.

There are several interesting future research directions. We assumed a deterministic demand in our model whereas in real world this cannot be a valid assumption. Moreover, we assumed an inspection and discard policy but in some industries like automotive and electronics industry this can be considered as inspection and fix policy where items defected after detecting can be repaired. Another good extension is to develop other metaheuristic techniques such as Genetic Algorithm or Tabu Search to compare their effectiveness with SA algorithm.

## Acknowledgements

This research was partially supported by NSF grant no. 08-068.

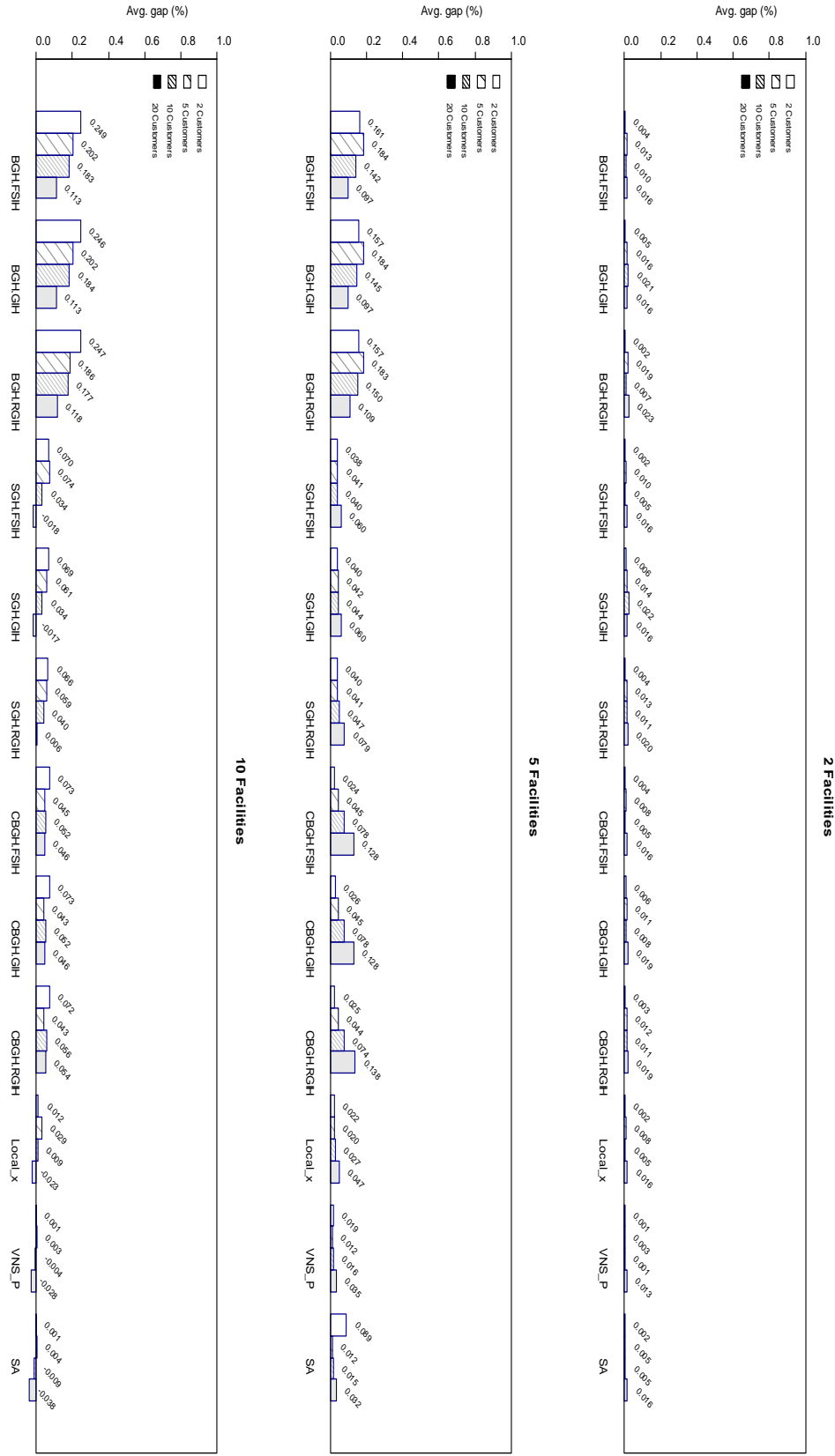

Figure 12 Gap of the algorithms under different number of facilities

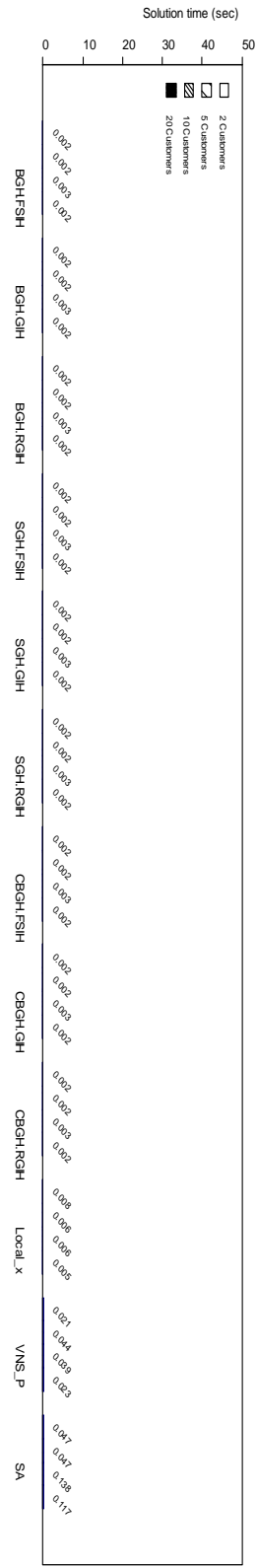
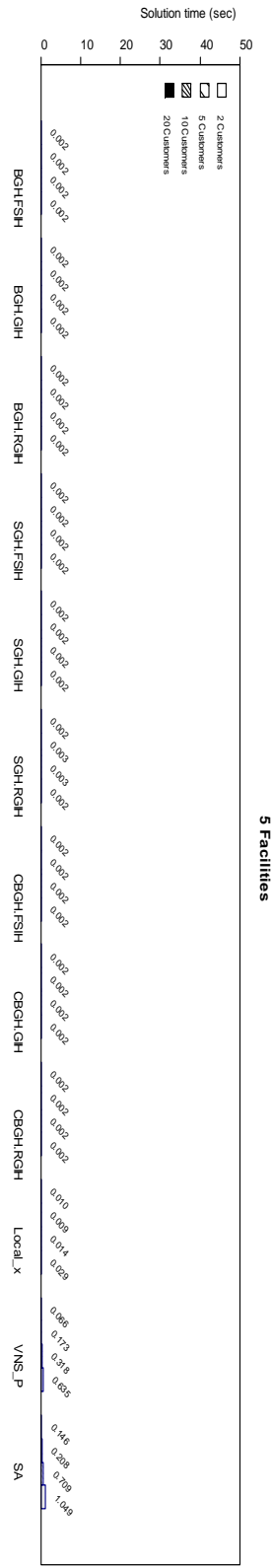
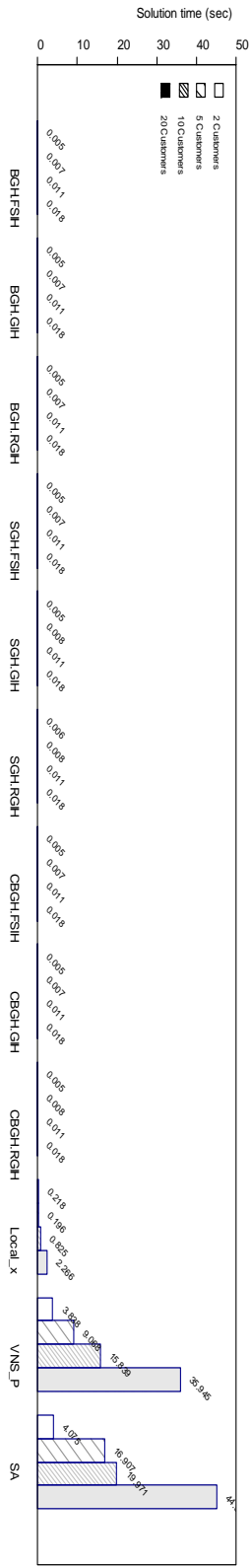

Figure 13 Solution time of the algorithms under different number of facilities